\documentclass{article}

\usepackage[russian,
english]{babel}

\usepackage{maa-monthly-my}
\final

\usepackage{hyperref}
\hypersetup{
colorlinks=true,citecolor=black,linkcolor=black,anchorcolor=black,filecolor=black,menucolor=black,urlcolor=black,
pdftitle={Moon in a puddle and the four-vertex theorem},
pdfauthor={Anton Petrunin and Sergio Zamora Barrera with artwork by Ana Cristina Chávez Cáliz}
}

\usepackage{caption}
\usepackage{wrapfig}
\usepackage[some,center]{background}
\backgroundsetup{
scale=.6,
opacity=1,
angle=0,
vshift=-40mm,
contents={%
  \includegraphics[width=\paperwidth
  ]{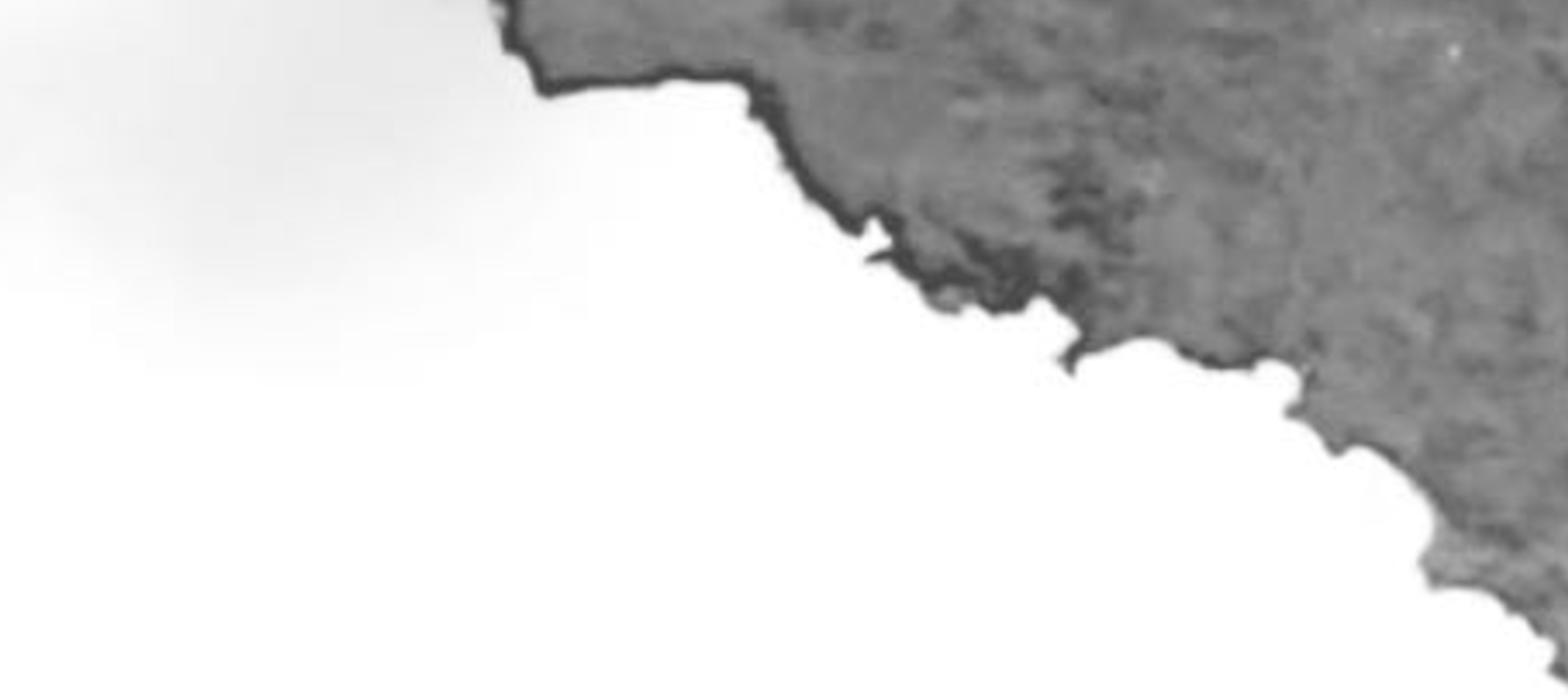}
  }%
}

\newcommand*{\z}[1]{#1\nobreak\discretionary{}%
            {\hbox{$\mathsurround=0pt #1$}}{}}

\theoremstyle{theorem}
\newtheorem{theorem}{Theorem}
\newtheorem*{keylemma}{Key lemma}
\newtheorem*{exercise}{Exercise}
\newtheorem*{advancedexercise}{Advanced exercise}
 
\theoremstyle{definition}

\begin{document}
\final
\title{Moon in a puddle and\\ the four-vertex theorem}
\markright{Moon in a puddle}
\author{Anton Petrunin and Sergio Zamora Barrera\\ with artwork by Ana Cristina Chávez Cáliz}

\maketitle
\BgThispage
\thispagestyle{empty}
\begin{abstract}
We present a proof of the moon in a puddle theorem, and use its key lemma to prove a generalization of the four-vertex theorem.
\end{abstract}

\newpage

\section*{INTRODUCTION.}

The theorem about the Moon in a puddle provides the simplest meaningful example of a local-to-global theorem which is mainly what differential geometry is about.
Yet, the theorem is surprisingly not well-known.
This paper aims to redress this omission by calling attention to the result and applying it to a well-known theorem

\section*{MOON IN A PUDDLE.}
The following question was initially asked by Abram Fet and solved by Vladimir Ionin and German Pestov \cite{pestov-ionin}.

\begin{theorem}\label{thm:moon-orginal}
Assume $\gamma$ is a simple closed smooth regular plane curve with curvature bounded in absolute value by~1.
Then the region surrounded by $\gamma$ contains a unit disc.
\end{theorem}

We present the proof from our textbook \cite{petrunin-zamora} which is 
a slight improvement of the original proof.
Both proofs work under the weaker assumption that the signed curvature is at most one, assuming that the sign is chosen suitably.
A more general statement for a barrier-type bound on the curvature was given by Anders Aamand, Mikkel Abrahamsen, and Mikkel Thorup~\cite{aamand-abrahamsen-thoru}.
There are other proofs. 
One is based on the curve-shortening flow; it is given by Konstantin Pankrashkin  \cite{pankrashkin}.
Another one uses cut locus; it is sketched by Victor Toponogov 
\cite[Problem 1.7.19]{toponogov}; see also \cite{petrunin-2020}.

Let us mention that an analogous statement for surfaces does not hold --- there is a solid body $V$ in the Euclidean space bounded by a smooth surface whose principal curvatures are bounded in absolute value by 1 such that $V$ does not contain a unit ball; moreover one can assume that $V$ is homeomorphic to the 3-ball.
Such an example can be obtained by inflating a nontrivial contractible 2-complex in $\mathbb{R}^3$ 
(Bing's house constructed in \cite{bing} would do the job).
This problem is discussed by Abram Fet and Vladimir Lagunov \cite{lagunov-2,lagunov-fet}; see also \cite{petrunin-zamora}.

\medskip

A path $\gamma\colon [0,1]\z\to \mathbb{R}^2$ such that $\gamma (0) = \gamma (1)$ will be called a \emph{loop};
the point $\gamma (0)$ is called the \emph{base} of the loop.
A loop is \emph{smooth}, \emph{regular}, and \emph{simple} if it is smooth and regular in $[0,1]$, and injective in the open interval $(0,1)$.

Let us use the term \emph{circline} as a shorthand for a \emph{circle or line}.
Note that the osculating circline of a smooth regular curve is defined at each of its points --- there is no need to assume that the curvature does not vanish.

Suppose that $\gamma$ is a closed simple smooth plane loop.
We say that a circline $\sigma$ \emph{supports} $\gamma$ at a point $p$ if the point $p$ lies on both $\sigma$ and $\gamma$, and the cicrline $\sigma$ lies in one of the closed regions that $\gamma$ cuts from the plane.
If furthermore this region is bounded, then  we say that $\sigma$ \emph{supports} $\gamma$ \emph{from inside}.
Otherwise, we say that $\sigma$ \emph{supports} $\gamma$ \emph{from the outside}.

\begin{keylemma}\label{thm:moon}
Assume $\gamma$ is a simple smooth regular plane loop.
Then at one point of $\gamma$ (distinct from its base), its osculating circle $\sigma$ supports $\gamma$ from inside.
\end{keylemma}

Spherical and hyperbolic versions of this lemma were given in \cite[Lemma 8.2]{panov-petrunin} and \cite[Proposition 7.1]{alexakis-mazzeo} respectively.

\medskip\noindent\textit{Proof of the theorem modulo the key lemma.}
Since $\gamma$ has absolute curvature of at most~1, each osculating circle has radius of at least 1.

According to the key lemma, one of the osculating circles $\sigma$ supports $\gamma$ from inside.
In this case, $\sigma$ lies inside $\gamma$, whence the result.
\qed

\medskip\noindent\textit{Proof of the key lemma.}
Denote by $F$ the closed region surrounded by $\gamma$.
Arguing by contradiction,
assume that the osculating circle at each point $p\ne p_0$ on $\gamma$ does not lie in~$F$.
Given such a point $p$, let us consider the maximal circle $\sigma$ that lies entirely in $F$ and is tangent to $\gamma$ at $p$.
The circle $\sigma$ will be called the {}\emph{incircle} of $F$ at $p$.
\begin{figure}[!ht]
\vskip-3mm
\centering
\includegraphics{pic-32}
\vskip0mm
\end{figure}

Note that the curvature of the incircle $\sigma$ has to be strictly larger than the curvature of $\gamma$ at $p$, hence there is a neighborhood of $p$ in $\gamma$ that intersects $\sigma$ only at $p$.
Further note that the circle $\sigma$ has to touch $\gamma$ at another point at least;
otherwise, we could increase $\sigma$ slightly while keeping it inside $F$.

Choose a point $p_1\ne p_0$ on $\gamma$, and  let $\sigma_1$ be the incircle at $p_1$.
Choose an arc $\gamma_1$ of $\gamma$ from $p_1$ to a first point $q_1$ on $\sigma_1$.
Denote by $\hat\sigma_1$ and $\check\sigma_1$ the two arcs of $\sigma_1$ from $p_1$ to $q_1$ such that the cyclic concatenation of $\hat\sigma_1$ and $\gamma_1$ surrounds~$\check\sigma_1$.

Let $p_2$ be the midpoint of $\gamma_1$, and $\sigma_2$ be the incircle at $p_2$.

Note that $\sigma_2$ cannot intersect $\hat\sigma_1$.
Otherwise, if $s$ is a point of the intersection, then $\sigma_2$ must have two more common points with $\check\sigma_1$, say $x$ and $y$ --- one for each arc of $\sigma_2$ from $p_2$ to $s$.
Therefore $\sigma_1\z=\sigma_2$ since these two circles have three common points: $s$, $x$, and $y$. 
On the other hand, by construction, $p_2\in \sigma_2$ and $p_2\notin \sigma_1$ --- a contradiction.

\begin{wrapfigure}{r}{37 mm}
\vskip-2mm
\centering
\includegraphics{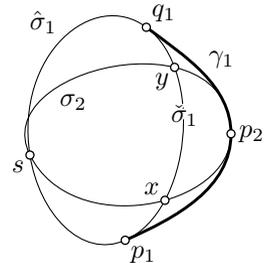}
\caption*{Two ovals pretend to be circles.}
\vskip0mm
\end{wrapfigure}

Recall that $\sigma_2$ has to touch $\gamma$ at another point.
From above it follows that it cannot touch $\gamma \setminus \gamma_1$, and therefore we can choose an arc $\gamma_2$ in $\gamma_1$ that runs from $p_2$ to a first point $q_2$ on $\sigma_2$.
Since $p_2$ is the midpoint of $\gamma_1$, we have that
\[\mathop{\rm length}\nolimits \gamma_2< \tfrac12\cdot\mathop{\rm length}\nolimits\gamma_1.\leqno({*})\]

Repeating this construction recursively,
we obtain an infinite sequence of arcs $\gamma_1\supset \gamma_2\supset\dots$;
by $({*})$, we also get that 
\[\mathop{\rm length}\nolimits\gamma_n\to0\quad\text{as}\quad n\to\infty.\] 
Therefore the intersection $\gamma_1\cap\gamma_2\cap\dots$
contains a single point; denote it by $p_\infty$.

Let $\sigma_\infty$ be the incircle at $p_\infty$; it has to touch $\gamma$ at another point, say $q_\infty$.
The same argument as above shows that $q_\infty\in\gamma_n$ for any $n$.
It follows that $q_\infty =p_\infty$ --- a contradiction.
\qed

\begin{exercise}\label{ex:moon-rad}
Assume that a closed smooth regular curve (possibly with self-intersections) $\gamma$  lies in a figure $F$ bounded by a closed simple plane curve.
Suppose that $R$ is the maximal radius of a disc contained in $F$.
Show that the absolute curvature of $\gamma$ is at least $\tfrac1R$ at some parameter value.
\end{exercise}

\section*{FOUR-VERTEX THEOREM.}
Recall that a \emph{vertex} of a smooth regular curve is defined as a critical point of its signed curvature;
in particular, any local minimum (or maximum) of the signed curvature is a vertex.
For example, every point of a circle is a vertex.

The classical four-vertex theorem says that \emph{any closed smooth regular plane curve without self-intersections has at least four vertices}.
It has many different proofs and generalizations.
A very transparent proof was given by Robert Osserman \cite{osserman}; his paper contains a short account of the history of the theorem.

Note that if an osculating circline $\sigma$ at a point $p$ supports $\gamma$, then $p$ is a vertex.
The latter can be checked by direct computation, but it also follows from the Tait--Kneser spiral theorem \cite{ghys-tabachnikov-timorin}.
It states that the \emph{osculating circlines of a curve with monotonic curvature are disjoint and nested};
in particular, none of these circlines can support the curve.
Therefore the following theorem is indeed a generalization of the four-vertex theorem:

{

\begin{wrapfigure}{r}{33 mm}
\vskip-4mm
\centering
\includegraphics{pic-63}
\vskip0mm
\end{wrapfigure}

\begin{theorem}\label{thm:4-vert}
Any smooth regular simple plane curve is supported by its osculating circlines at 4 distinct points; two from inside and two from outside.
\end{theorem}

\medskip\noindent\textit{Proof.}
According to the key lemma, there is a point $p\in\gamma$ such that its osculating circle supports $\gamma$ from inside.

The curve $\gamma$ can be considered as a loop with $p$ as its base.
Therefore the key lemma implies the existence of another point $q$ with the same property.

This shows the existence of two osculating circles that support $\gamma$ from inside;
it remains to show the existence of two osculating circles that support $\gamma$ from outside.

Let us apply to $\gamma$ an inversion with respect to a circle whose center lies inside~$\gamma$.
Then the obtained curve $\gamma_1$ also has  two osculating circles that support $\gamma_1$ from inside.

}

Note that these osculating circlines are inverses of the osculating circlines of $\gamma$.
Indeed, the osculating circline at a point $x$ can be defined as the unique circline that has second order of contact with $\gamma$ at $x$.
It remains to note that inversion, being a local diffeomorphism away from the center of inversion, does not change the order of contact between curves.

Note that the region lying inside $\gamma$ is mapped to the region outside $\gamma_1$ and the other way around.
Therefore these two new circlines correspond to the osculating circlines supporting $\gamma$ from outside.
\qed

\begin{advancedexercise}\label{ex:curve-crosses-circle}
Suppose $\gamma$ is a closed simple smooth regular plane curve and $\sigma$ is a circle.
Assume $\gamma$ crosses $\sigma$ at the points $p_1,\dots,p_{2{\cdot}n}$ and these points appear in the same cycle order on $\gamma$ and on $\sigma$.
Show that $\gamma$ has at least $2{\cdot}n$ vertices.
\end{advancedexercise}

\begin{figure}[!ht]
\begin{minipage}{.48\textwidth}
\centering
\includegraphics{pic-65}
\end{minipage}\hfill
\begin{minipage}{.48\textwidth}
\centering
\includegraphics{pic-305}
\end{minipage}
\end{figure}

The order of the intersection points is important. 
An example with only 4 vertices and arbitrarily many intersection points can be guessed from the diagram on the right.

\begin{acknowledgment}{Acknowledgments.}
We wish to thank anonymous referees for thoughtful reading and insightful suggestions.
This work was supported by the National Science Foundation under Grant DMS-2005279; Simons Foundation under Grant \#584781.
\end{acknowledgment}

\end{document}